\documentclass[preprint,12pt]{article}

 \usepackage{graphics}
 \usepackage{graphicx}
 \usepackage{epsfig}
 \usepackage[all]{xy}

\usepackage{amssymb}
\usepackage{amsthm,amsmath}

\usepackage{latexsym, epsfig, psfrag,eepic,colordvi,bm}
\usepackage{graphics,color}
\usepackage{amsmath,amsfonts,amssymb,amscd}
\usepackage[all]{xy}
\usepackage{epstopdf}
\usepackage{cite}

\usepackage{amsthm,amsmath,amssymb}

\usepackage{graphicx}

\theoremstyle{plain}
\newtheorem{theorem}{Theorem}[section]
\newtheorem{lemma}[theorem]{Lemma}
\newtheorem{corollary}[theorem]{Corollary}
\newtheorem{proposition}[theorem]{Proposition}

\theoremstyle{definition}
\newtheorem{definition}[theorem]{Definition}
\newtheorem{example}[theorem]{Example}

\theoremstyle{remark}

\newcommand{\bai}{\hspace{6pt}}

\begin{document}

		

\title{A versatile combinatorial approach of studying products of long cycles in symmetric groups }

\author{Ricky X. F. Chen\\
	\small Biocomplexity Institute and Initiative, University of Virginia\\[-0.8ex]
	\small Charlottesville, VA 22901, USA\\
	\small\tt chenshu731@sina.com
}

\date{}
\maketitle


\begin{abstract}

In symmetric groups, studies of permutation factorizations or triples of permutations satisfying certain conditions have a long history. 
One particular interesting case is when two of the involved permutations are long cycles,
for which many surprisingly simple formulas have been obtained.
Here we
combinatorially enumerate the pairs of long cycles whose product has a given cycle-type and separates certain elements,
extending several lines of studies,
and we obtain general quantitative relations.
As consequences, in a unified way, we recover a number of results expecting simple combinatorial
proofs, including results of Boccara (1980),
Zagier (1995), Stanley (2011), F\'{e}ray and Vassilieva (2012), as well as Hultman (2014).

We obtain a number of new results as well. In particular, for the first time,
given a partition of a set, we obtain an explicit formula for the number of pairs of long 
cycles on the set such that the product of the long cycles does not mix the elements from distinct blocks of the partition and
has an independently prescribed number of cycles for each block of elements.
As applications, 
we obtain new explicit formulas concerning factorizations of any even permutation into long cycles and the first nontrivial explicit formula for computing strong separation probabilities
solving an open problem of Stanley (2010).

  

\end{abstract}



\section{Introduction}

Let $\mathfrak{S}_n$ denote the symmetric group on $[n]=\{1,2,\ldots, n\}$,
and let $C(\pi)$ denote the number of disjoint cycles of $\pi \in \mathfrak{S}_n$.
The set consisting of the lengths of these disjoint cycles is called the cycle-type of $\pi$. 
We can encode this set as an integer partition of $n$.
An integer partition $\lambda$ of $n$, denoted by $\lambda \vdash n$,
can be represented by a non-increasing integer sequence $\lambda=\lambda_1
\lambda_2\cdots$, where $\lambda_1 \geq \lambda_2 \geq \cdots, \; \sum_i \lambda_i=n$, or as $1^{m_1}2^{m_2}\cdots n^{m_n}$, where
we have $m_i(\lambda)$ of part $i$ and $\sum_i i m_i =n$.
 A cycle of length $k$ is called a $k$-cycle.
 We sometimes call an $n$-cycle on $[n]$ a long cycle. 
In addition, we denote the number of permutations of cycle-type $\lambda$ by $z_{\lambda}$.
It is well known that if $\lambda=1^{m_1} 2^{m_2} \cdots n^{m_n}$, then 
$$
z_{\lambda}=\frac{n!}{\prod_i i^{m_i}  m_i !}.
$$
We also denote the length of $\lambda$, i.e.,~the number of positive parts in $\lambda$, by $\ell(\lambda)$.

Factorizations of permutations or triples of permutations
satisfying certain conditions
have been extensively studied in different contexts.
A particular important and interesting case is when one of 
the involved permutations is a long cycle, e.g.,~enumeration of maps by
Walsh and Lehman~\cite{wl}, factorizations of permutations
and/or maps and/or genome rearrangement problems by
Boccara~\cite{boccara}, Stanley~\cite{stan}, 
Jackson~\cite{jack}, Zagier~\cite{zag}, Goupil and Schaeffer~\cite{gs}, Chapuy~\cite{chapuy}, Chapuy et al.~\cite{cff}, Bernardi~\cite{bernardi}, and Chen and Reidys~\cite{chr-1}, 
graph embeddings by Gross et al.~\cite{grt}, as well as studying the Euler characteristic of the moduli spaces of algebraic curves
by Harer and Zagier~\cite{hz}. Most of the related results in the field
 rely either partially or totally on character theoretic approaches (e.g.,~\cite{stan,jack,zag, gs})
 or integral approaches (e.g.,~\cite{boccara, hz})
 in earlier decades, while progresses from combinatorial approaches have been
 made very recently (e.g.,~\cite{CV, gs,chapuy, cff, bernardi, chr-1}). 
It is generally hard to obtain explicit and simple counting formulas.
However, if two of the involved permutations are long cycles~\cite{bdms,boccara,bona-pittel,berwei,can,cms,chen3,chen4,fr,fv,stan, stan3, walkup},
we may have very nice formulas for most of studied problems, at
least much simpler than those of the general case.
For instance, we have the following results that can be
clearly stated without requiring additional notation and definitions.
(We shall combinatorially prove them all, providing either the first combinatorial proof or
probably the most simple one.)

\begin{theorem}[Zagier~\cite{zag}, Stanley~\cite{stan3}]\label{thm:1t1}
The number of $n$-cycles $s$ such that the product $(1\bai 2\bai \cdots \bai n)\, s$ has $k$ cycles is $\frac{2}{n(n+1)} C(n+1,k)$,
where $C(n,k)$ stands for the signless Stirling number of the first kind, i.e., the number of permutations on $[n]$ with $k$ cycles.

\end{theorem}

\begin{theorem}[Hultman~\cite{hultman}]\label{thm:1t2}
The expected number of $k$-cycles
in the product of two random long cycles on $[n]$ is $\frac{(-1)^{k+1}}{k {n-1 \choose k}} + \frac{1}{k}$.
\end{theorem}

\begin{theorem}[Boccara~\cite{boccara}]\label{thm:1t3}
The number of different factorizations of a fixed even permutation on $[n]$
of cycle-type $k^1 (n-k)^1$ into two $n$-cycles is given by $\frac{2(n-1)!}{n+1} \Bigl( 1-\frac{(-1)^k}{{n \choose k}}\Bigr)$.
\end{theorem}

\begin{theorem}\label{thm:1t4}
The number of pairs
of long cycles whose product does not mix the elements in $[k]$
and the elements in $[n]\setminus [k]$ in any of its cycles is $k!(n-k)!(n-2)!$.
More generally, for an integer composition $\alpha=(\alpha_1, \ldots, \alpha_k)$
of $n$, the corresponding number
is 
$$
\frac{(n-1)! }{n+1-k} \prod_{t=1}^k \alpha_t ! \; .
$$
\end{theorem}

\begin{theorem}\label{thm:1t6}
The number of factorizations of a fixed even permutation of cycle-type $\lambda=\lambda_1 \lambda_2 \cdots \lambda_k \vdash n$
into two $n$-cycles
is given by
\begin{align*}
2(n-1)!  \sum_{j_2+\cdots +j_k =l, \atop l\geq 0, \; \lambda_t > j_t\geq 0}  (-1)^l  \frac{l!}{(\lambda_1+l+1)_{l+1}} {\lambda_2 \choose j_2} \cdots {\lambda_k \choose j_k}\; .
\end{align*}
\end{theorem}

In F\'{e}ray and Vassilieva~\cite{fv}, a refined problem of the one considered in Theorem~\ref{thm:1t1} was also considered, that is,
enumerating the pairs of long cycles whose product has a given cycle-type.
A simple relation of these refined numbers was obtained in~\cite{fv} by 
counting some colored permutations first and then by some algebraic computations in the ring of symmetric functions.

Separation probabilities for products of permutations
 were studied in Bernardi, Du, Morales, and Stanley~\cite{bdms}, where
a special case is concerned with the probability of having the elements in $[m]$
contained in distinct cycles of the product of two uniformly randomly chosen $n$-cycles.
It was proved~\cite{bdms, stan2} that the separation probability is given by
\begin{align*}
\begin{cases}
\frac{1}{m!} , & \mbox{ if $n-m$ is odd},\\
\frac{1}{m!}+\frac{2}{(m-2)! (n+1-m)(n+m)} , & \mbox{if $n-m$ is even.}
\end{cases}
\end{align*}
A more general case of separation probabilities studied was about studying the
pairs of long cycles such that each cycle of their product may
contain elements from a particular subset of $[n]$. (See a more precise description later.)

Although these formulas concerning products of two long cycles are simple,
simple combinatorial proofs are not necessarily immediately available.
Recently, the author obtained some analogues of the above Zagier-Stanley result in
the context of studying separation probabilities~\cite{chen3}.
For example, it was proved that the number of pairs of $n$-cycles whose product has $k$ cycles and separates the elements in $[m]$ is given by 
$$
\frac{2 (n-1)! }{(n+m)(n+1-m)} C_m(n+1,k),
$$
where $C_m(n,k)$ is the number of permutations on $[n]$
with $k$ cycles and the elements in $[m]$ separated, i.e.,~an analogue of $C(n,k)$.
Based on the analogue, the aforementioned formulas of separation probabilities immediately follow.
In addition, in~\cite{chen4}, with simpler combinatorial arguments, the author also obtained the relation previously obtained by
F\'{e}ray and Vassilieva~\cite{fv}. Both works were based on extending the plane permutation framework that was first introduced in Chen and Reidys~\cite{chr-1}
in order for studying one-face hypermaps as well as genome rearrangement problems.
Accordingly, it was speculated that the approach based on the plane permutation framework may take us even further,
which motivated the present work.
For instance, we shall prove the above theorems in a unified way later.

The main problem studied in this paper is described as follows.
Let $k>0$ and $\alpha=(\alpha_1, \alpha_2, \ldots, \alpha_k)$ be an integer composition of $n$, i.e.,~$\sum_i \alpha_i =n$ and $\alpha_i>0$. We write $\alpha \models n$.
Let $B_i=\{\sum_{j=0}^{i-1} \alpha_j+1, \ldots, \sum_{j=0}^{i-1} \alpha_j+\alpha_i\} \subseteq [n]$ where we assume $\alpha_0=0$.
A permutation $\pi$ on $[n]$ is called $\alpha$-separated if the elements in every cycle of $\pi$
are coming from the same $B_i$ for some $1\leq i \leq k$.
We shall be interested in the pairs of $n$-cycles whose product
is $\alpha$-separated, and we shall refer to these pairs as $\alpha$-separating pairs.

This problem appears to be a special case of the separation probabilities considered in~\cite{bdms},
where any $\alpha \models m \leq n$ was considered, i.e.,~there are some
free elements. However, we shall discuss later that it suffices to merely study the case $\alpha \models n$.
We shall actually begin with studying $\alpha$-separating pairs whose product
has a given cycle-type. Therefore,
our studied problem here can be viewed as extending a number of previously cited works at
the same time.

\section{Review of plane permutations}\label{sec2}
As a new way of representing one-face hypermaps, plane permutations were introduced in~\cite{chr-1} for the first time. As results, several recurrences counting one-face hypermaps or factorizations of a long cycle were obtained, and a combinatorial proof of the Zagier-Stanley result was given, answering a call of Stanley~\cite{stan3}
for a combinatorial proof. Plane permutations also provided a unified simple framework studying transposition and block-interchange distances of permutations, as well as reversal distances of signed permutations. For the latter application to the genome rearrangement problems, a general lower bound was obtained, for which the well-known Bafna-Pevzner's lower bound~\cite{trans} and Christie's formula~\cite{christie} based on cycle graphs are equivalent to evaluations at a special point. Plane permutations were later employed to study the local genus
distribution problem of graph embeddings, where the local genus distribution
of reembedding a single vertex of a graph was fully solved and the local genus 
distribution was shown to be indeed log-concave~\cite{chr-2}.				
Let us review some notation and results in this section.

\begin{definition}
	A \emph{plane permutation} on $[n]$ is a pair $\mathfrak{p}=(s,\pi)$ where $s=(s_i)_{i=0}^{n-1}$
	is an $n$-cycle and $\pi$ is an arbitrary permutation on $[n]$.
	Given $s=(s_0~s_1~\cdots ~s_{n-1})$,
a plane permutation $\mathfrak{p}=(s,\pi)$ is represented by a two-row array:
\begin{equation}
\mathfrak{p}=\left(\begin{array}{ccccc}
s_0&s_1&\cdots &s_{n-2}&s_{n-1}\\
\pi(s_0)&\pi(s_1)&\cdots &\pi(s_{n-2}) &\pi(s_{n-1})
\end{array}\right).
\end{equation}
The permutation $D_{\mathfrak{p}}$ induced by the diagonal-pairs (cyclically), i.e.,~$D_{\mathfrak{p}}(\pi(s_{i-1}))=s_i$ for $0<i< n$, and
$D_{\mathfrak{p}}(\pi(s_{n-1}))=s_0$, is called the \emph{diagonal} of $\mathfrak{p}$.

\end{definition}\label{2def1}

We sometime refer to $s,\, \pi, \, D_{\mathfrak{p}}$ respectively as the upper horizontal, the vertical and the diagonal. The following lemma is obvious but important.

\begin{lemma}\label{observation}
$D_{\mathfrak{p}}=
s \pi^{-1}$.
\end{lemma}
Based on Lemma~\ref{observation}, we can see that studying factorizations of 
permutations is equivalent to studying plane permutations.

In a permutation $\pi$ on $[n]$, $i$ is called an \emph{exceedance} if $i<\pi(i)$ following the natural order and an \emph{anti-exceedance} otherwise.
Exceedances and anti-exceedances are among the most well-known 
permutation statistics.
Note that $s$ induces a linear order $<_s$,
where $a<_s b$ if $a$ appears before $b$ in $s$ from left to right (with the left most element $s_0$).
Without loss of generality, we always assume $s_0=1$ unless explicitly stated otherwise.
These concepts then can be generalized for plane permutations as follows:
\begin{definition}\label{2def2}
	For a plane permutation $\mathfrak{p}=(s,\pi)$, an element $s_i$ is called an
	\emph{exceedance} of $\mathfrak{p}$ if $s_i<_s \pi(s_i)$, and an \emph{anti-exceedance} if $s_i\ge_s \pi(s_i)$.
\end{definition}

In the following,
any comparison of elements in $s,~\pi$ and $D_{\mathfrak{p}}$ references the order $<_s$.
Obviously, each $\mathfrak{p}$-cycle contains at least one anti-exceedance as it contains
a minimum, $s_i$, for which $\pi^{-1}(s_i)$ is an anti-exceedance. We call these trivial anti-exceedances
and refer to a \emph{non-trivial anti-exceedance} as an NTAE. Furthermore, in any cycle of length
greater than one, its minimum is always an exceedance.

\begin{example}
	For the plane permutation 
	\begin{equation*}
	\mathfrak{p}=\left(\begin{array}{cccccc}
	1& 5& 4 &6 & 2 & 3\\
	5& 4 &1&3 & 6 & 2
	\end{array}\right),
	\end{equation*}
	$1$ is an exceedance, $2$ is an anti-exceedance and also an NTAE.
\end{example}

Let $\mathfrak{p}=(s,\pi)$ be a plane permutation. 
A \emph{diagonal block} of $\mathfrak{p}$ is a set of consecutive diagonal-pairs.
A \emph{transposition action} on the diagonal of $\mathfrak{p}$
transposes two adjacent diagonal blocks of $\mathfrak{p}$.
Specifically,
for a sequence $h=(i,j,k)$ such that $i\leq j<k$
and $\{i,j,k\}\subset [n-1]$, if we transpose the two diagonal-blocks determined by the continuous segments $[s_i,s_j]$ and $[s_{j+1},s_k]$, we obtain a new two-row array $\mathfrak{p}^h=(s^h,\pi^h)$:
\begin{eqnarray*}
\left(
\vcenter{\xymatrix@C=0pc@R=1pc{
\cdots & s_{i-1}  & s_{j+1}\ar@{--}[dl] &\cdots & s_{k-1}&\bai s_k\ar@{--}[dl] \bai\bai& s_i\ar@{-}[dl] &\cdots & s_{j-1} & s_{j}\ar@{-}[dl] & s_{k+1}  &\cdots\\
\cdots  & \pi(s_{j}) & \pi(s_{j+1}) & \cdots & \pi(s_{k-1})& \pi(s_{i-1}) & \pi(s_i) & \cdots\hspace{-0.5ex} & \pi(s_{j-1}) & \pi(s_{k})  & \pi(s_{k+1})&\cdots
}}
\right).
\end{eqnarray*}
Comparing $\mathfrak{p}$ and $\mathfrak{p}^h$, we have the following observations:
\begin{itemize}
	\item they have the same diagonal;
	\item the upper horizontals $s$ and $s^h$ differ by a transposition of the two continuous segments $[s_i,s_j]$ and $[s_{j+1},s_k]$;
	\item the maps $\pi$ and $\pi^h$ only differ at the images of the elements $s_{i-1}$, $s_j$, and $s_k$.
\end{itemize}
Thus, the transposition actions on the diagonal provide
a natural viewpoint on how different factorizations of the diagonal
into a long cycle (the upper horizontal) and another permutation (the vertical)
relate to each other. 
In particular,
the above last bullet implies that all components other than those containing the mentioned three elements of $\pi$ will be completely carried over to $\pi^h$ without any changes.
For those components containing the three elements, the three elements serve as certain breakpoints, where
the induced segments will be re-pasted in a certain way, depending on the distribution of the elements $s_{i-1}$, $s_j$, and $s_k$ in the components of $\pi$. 
Note that $\pi$ and $\pi^h$ must have the same parity.
Thus, the difference of the number of cycles in $\pi^h$ and $\pi$ is contained in $\{2, 0,-2\}$.
The NTAEs of $\mathfrak{p}$ can help us
to identify the transposition actions which change the number of cycles
by exactly two.

Let $\mu,~\lambda$ be two integer partitions of $n$.
We denote $\mu\rhd_{k} \lambda$ if $\mu$ can be obtained from $\lambda$ by splitting one part into $k$ parts,
or equivalently, $\lambda$ from $\mu$ by merging $k$ parts into one part.
Let $\kappa_{\mu,\lambda}$ be the number of different ways of merging $k$ parts of $\mu$
in order to obtain $\lambda$ provided that $\mu \rhd_k \lambda$.
Note that we differentiate two parts of $\mu$ even if the two parts are of the same value.
For example, for $\mu=1^2 2^2$ and $\lambda=1^1 2^1 3^1$, we have
$\kappa_{\mu, \lambda}=4$.

Let $D$ be a fixed permutation on $[n]$.
We consider the set of plane permutations
$\mathfrak{p}=(s, \pi)$ where the diagonal is $D$ and the vertical
is of cycle-type $\lambda$, or equivalently the set of factorizations of $D$ into 
a long cycle $s$ on $[n]$ and a permutation $\pi^{-1}$ on $[n]$ of cycle-type
$\lambda$, i.e.,~$D=s\pi^{-1}$.
Denote this set by $\tilde{U}_{\lambda}^D$. 
By studying the transposition actions on the diagonal of plane permutations
with the assistance of NTAEs,
the following result has been obtained.

 \begin{proposition}[Chen\&Reidys~\cite{chr-1}]\label{prop:general}
 Let $\tilde{Y}_1$ be the set of pairs $(\mathfrak{p},\epsilon)$ where $\mathfrak{p}\in \tilde{U}_{\lambda}^D$ and $\epsilon$
 is an NTAE of $\mathfrak{p}$. Let $\tilde{Y}_2$ be the set of plane permutations $\mathfrak{p}\in \bigcup_{j\geq 1, \; \mu\rhd_{2j+1} \lambda}  \tilde{U}_{\mu}^D$ where there are $2j+1$ marked cycles in $\mathfrak{p}$ if $\mathfrak{p}\in  \tilde{U}_{\mu}^D$ and $\mu\rhd_{2j+1} \lambda$
 such that when treating the $2j+1$ marked cycles as a single cycle the cycle-type of the vertical of $\mathfrak{p}$ is $\lambda$. 
 Then there is a bijection between $\tilde{Y}_1$ and $\tilde{Y}_2$.
 \end{proposition}
 
 The main idea behind the above bijection between $\tilde{Y}_1$ and $\tilde{Y}_2$ can be briefly summarized here.
 From a given pair $(\mathfrak{p},\epsilon)$ in $\tilde{Y}_1$, the NTAE $\epsilon$ determines a transposition on the diagonal 
 of $\mathfrak{p}$ such that the vertical of the resulting plane permutation after the transposition is
 obtained by splitting the cycle containing $\epsilon$ of $\mathfrak{p}$ into three cycles.
Obviously, $\epsilon$ is still contained in one of the three cycles.
Depending on whether $\epsilon$ is still an NTAE of the resulting plane permutations, additional
transpositions can be applied until $\epsilon$ is not an NTAE anymore.
Eventually, the original cycle containing $\epsilon$ will split into $2j+1$ cycles for some $j>0$ which will be marked
and the resulting plane permutation has a cycle-type $\mu \rhd_{2j+1} \lambda$ in the vertical.
Conversely, from a given element in $\tilde{Y}_2$, there is a unique way to merge the marked cycles
into one single cycle and create an NTAE.
We refer to~\cite{chr-1} for details.

The above bijection was motivated by the vertex slicing/gluing bijection on one-face maps of Chapuy~\cite{chapuy}.
 However, once we had the two-row array formulation of plane permutations, 
 it was in fact the natural transposition action on the diagonal, or an even broader perspective, rearrangement of the diagonal-pairs, that were first studied, due to their clear potential applications to the transposition, block-interchange and reversal distances of genome sequences.
 It turned out that the slicing/gluing operations are hiding there as two particular cases among others (see~\cite[Lemma~$7$]{chr-1}),
 somehow resolving the mystery of the slicing/gluing bijection~\cite{chapuy}.

 It should be easy to see that the cardinality $|\tilde{U}_{\lambda}^D|$ only depends
 on the cycle-type $\eta$ of $D$.
 Let $\tilde{U}_{\lambda}^{\eta}$ denote the set of plane permutations on $[n]$
where the diagonal is of cycle-type $\eta$ and the vertical is of cycle-type $\lambda$. 
Then, we have $|\tilde{U}_{\lambda}^{\eta}|= z_{\eta} |\tilde{U}_{\lambda}^D|$.
We always assume $\ell(\lambda)+\ell(\eta)$ has the same parity as $n+1$.
Otherwise we know $\tilde{U}_{\lambda}^{\eta}=\varnothing$.
We denote $p^{\eta}_{\lambda}= |\tilde{U}_{\lambda}^{\eta}|$,
in particular, for $\eta=n^1$, we define $p^{(n)}_{\lambda}:=p^{\eta}_{\lambda}$.
Let $Ne(\mathfrak{p})$ denote the number of NTAEs of $\mathfrak{p}$.
Note that for $\mathfrak{p} \in \tilde{U}^{\eta}_{\lambda}$ having $a$ exceedances,
we have $Ne(\mathfrak{p})=n-\ell(\lambda)-a$.
Then, as a consequence of the bijection in Proposition~\ref{prop:general},
we obtain the quantitative relations below.

 \begin{proposition}[Chen\&Reidys~\cite{chr-1}]
Let $\lambda, \eta \vdash n$ and $p^{\eta}_{\lambda,a}$ be the number of $\mathfrak{p} \in \tilde{U}^{\eta}_{\lambda}$ such that $\mathfrak{p}$
has $a$ exceedances. Then
we have
\begin{align}
\sum_{a\geq 0} \bigl(n-\ell(\lambda)-a \bigr) p^{\eta}_{\lambda,a} &=\sum_{ \mu \rhd_{2i+1} \lambda, \; i>0} \kappa_{\mu,\lambda} p^{\eta}_{\mu} \label{eq:gen}\; ,\\
\sum_{a\geq 0} \bigl(n-\ell(\eta)-a \bigr) p^{\lambda}_{\eta, a} &=\sum_{ \mu \rhd_{2i+1} \eta, \; i>0} \kappa_{\mu,\eta} p^{\lambda}_{\mu} \label{eq:gen-reflection}\; .
\end{align}
\end{proposition}

We remark that these equations in the above corollary are inherently 
filtered out (or avoided) if following the map (and bicolored map~\cite{chapuy,walsh2}) perspective.
Because on the one hand, all plane permutations corresponding to maps (i.e.,~the diagonal being a fixed point free involution) have the same fixed number
(roughly speaking, $2g$ for $g$ being the genus)
of NTAEs, such equations never appear in the first place; On
the other hand, these equations do not really provide `valid'
recurrences from an enumeration perspective. (In order for obtaining valid recurrences, we have to apply a
sort of `reflection principle' to clear the parameter `a', as will be shown shortly.)
However, these `invalid' recurrences are actually 
the most valuable ingredients to make our approach productive.

Note that there is a one-to-one correspondence between $\tilde{U}^{\eta}_{\lambda}$
and $\tilde{U}^{\lambda}_{\eta}$, as if $\mathfrak{p}=(s, \pi) \in \tilde{U}^{\eta}_{\lambda}$,
then $\mathfrak{p'}=(s^{-1}, D_{\mathfrak{p}}^{-1}) \in \tilde{U}^{\lambda}_{\eta}$.
There is a also nice relation stated in the following lemma.

 \begin{lemma}[Chen\&Reidys~\cite{chr-1}]\label{reflection}
Let $\mathfrak{p}=(s,\pi)$ be a plane permutation with diagonal $D_{\mathfrak{p}}$,
and let $\mathfrak{p'}=(s^{-1}, D_{\mathfrak{p}}^{-1})$. Then,
\begin{align}\label{eq:reflection}
Ne(\mathfrak{p})+Ne(\mathfrak{p'})=n+1-C(\pi)-C(D_{\mathfrak{p}}) \; .
\end{align}
\end{lemma}

Now, based on eq.~\eqref{eq:gen}, eq.~\eqref{eq:gen-reflection}, and eq.~\eqref{eq:reflection},
we obtain
 \begin{proposition}[Chen\&Reidys~\cite{chr-1}]
Let $\lambda, \eta \vdash n$. Then,
we have
\begin{align}
\bigl( n+1-\ell(\lambda)-\ell(\eta)\bigr) p^{\eta}_{\lambda} &=\sum_{ \mu \rhd_{2i+1} \lambda, \; i>0} \kappa_{\mu,\lambda} p^{\eta}_{\mu} + \sum_{ \mu \rhd_{2i+1} \eta, \; i>0} \kappa_{\mu,\eta} p^{\mu}_{\lambda} \label{eq:gen-explicit}\; ,\\
\bigl( n+1-\ell(\lambda) \bigr) p^{(n)}_{\lambda} &=\sum_{\mu \rhd_{2i+1} \lambda, \; i>0} \kappa_{\mu,\lambda} p^{(n)}_{\mu} + (n-1)! z_{\lambda} \label{eq:long}\; .
\end{align}
\end{proposition}

In order for obtaining eq.~\eqref{eq:long} from eq.~\eqref{eq:gen-explicit}, we have used the fact that
any partition $\mu$ with $\ell(\mu)+\ell(\lambda)$ having the same parity as $1+\ell(\lambda)$ must satisfy $\mu\rhd_{2i+1} 1^n$ for some $i\geq 0$ and $\kappa_{\mu, 1^n}=1$,
and $\sum_{\mu \vdash n} p^{\lambda}_{\mu} =(n-1)! z_{\lambda}$.

 From the summarized bijection above, new bijections can be derived if there are some appropriate
 restrictions on the set of plane permutations under consideration. 
 The results in the rest of the paper are based on such a bijection.

\section{Separating pairs of long cycles}\label{sec3}

In this section, we present our general results based on a revised version
of the bijection described in Proposition~\ref{prop:general}. 				

We begin with introducing the following notation.
For $\alpha=(\alpha_1, \alpha_2, \ldots, \alpha_k) \models n$, we write ${\Lambda} \vdash \alpha$ if ${\Lambda}=(\lambda^{[1]}, \lambda^{[2]}, \ldots, \lambda^{[k]} )$ is a sequence of
integer partitions where $\lambda^{[i]} \vdash \alpha_i$.
We also define $z_{\Lambda}:=\prod_{i=1}^k z_{\lambda^{[i]}}$,
$\ell(\Lambda):=\sum_{i=1}^k \ell(\lambda^{[i]})$, and $m_i(\Lambda):= \sum_{j=1}^k m_i(\lambda^{[j]})$.
Moreover, for $\Upsilon \vdash \alpha$, we write $\Upsilon \rhd_{i,2j+1} \Lambda$ if 
$\Upsilon$ can be obtained from $\Lambda$ by splitting one part of $\lambda^{[i]}$ 
into $2j+1$ parts. The number $\kappa_{\Upsilon, \Lambda}$ is analogously defined.

An $\alpha$-separated permutation $\pi$ has
an $\alpha$-type $\Lambda$ if the elements in $B_i$ have a cycle-type $\lambda^{[i]}$.
Let $U_{\Lambda}^{\eta}$ denote the set of plane permutations on $[n]$
where the diagonal is of cycle-type $\eta$ and the vertical has an $\alpha$-type $\Lambda$. 
With the same reason, we always assume $\ell(\Lambda)+\ell(\eta)$ has the same parity as $n+1$
if not explicitly specified.

We denote $p^{\eta}_{\Lambda}= |U_{\Lambda}^{\eta}|$, and specially $p^{(n)}_{\Lambda}$ for $\eta=n^1$.
Note that the bijection described in the last section inherently preserves the property of the vertical 
being $\alpha$-separated in the forward direction.
As to the converse, we merely need to pay attention to
how we mark the cycles, i.e.,
we can only mark cycles containing elements from the same $B_i$.
 Hence, we can succinctly obtain the following proposition.

\begin{proposition}
Suppose $\alpha \models n,\; \Lambda \vdash \alpha$ and $\eta \vdash n$. 
Let $p^{\eta}_{\Lambda,a}$ be the number of $\mathfrak{p} \in U^{\eta}_{\Lambda}$ such that $\mathfrak{p}$
has $a$ exceedances. Then
we have
\begin{align}
\sum_{a\geq 0} \bigl(n-\ell(\Lambda)-a \bigr) p^{\eta}_{\Lambda,a} &=\sum_{i>0,\; j>0 \atop \Upsilon \rhd_{i, 2j+1} \Lambda} \kappa_{\Upsilon,\Lambda} p^{\eta}_{\Upsilon} \label{eq:gen-new}\; ,\\
\bigl( n+1-\ell(\Lambda)-\ell(\eta)\bigr) p^{\eta}_{\Lambda} &=\sum_{i>0, \; j>0 \atop \Upsilon \rhd_{i, 2j+1} \Lambda} \kappa_{\Upsilon,\Lambda} p^{\eta}_{\Upsilon}+ \sum_{ \mu \rhd_{2i+1} \eta, \; i>0} \kappa_{\mu,\eta} p^{\mu}_{\Lambda}  \label{eq:gen-explicit-new} \; ,\\
\bigl( n+1-\ell(\Lambda) \bigr) p^{(n)}_{\Lambda} &=\sum_{i>0, \; j>0 \atop \Upsilon \rhd_{i, 2j+1} \Lambda} \kappa_{\Upsilon,\Lambda} p^{(n)}_{\Upsilon} + (n-1)! z_{\Lambda} \label{eq:long-new}\; .
\end{align}
\end{proposition}

{\noindent\bf Remark.} The equations eq.~\eqref{eq:gen-new}, eq.~\eqref{eq:gen-explicit-new},
and eq.~\eqref{eq:long-new} respectively correspond to eq.~\eqref{eq:gen}, eq.~\eqref{eq:gen-explicit},
and eq.~\eqref{eq:long}. 
Although in the present paper, we focus on $\alpha$-separating
pairs of long cycles, eq.~\eqref{eq:gen-explicit-new} and its corresponding initial values allow us to enumerate $\alpha$-separating pairs w.r.t.~any arbitrary $\eta$. Explicit formulas for the corresponding initial values can be obtained
in the similar manner as the general case concerning separation probabilities discussed in~\cite{chen3} by enumerating certain labelled plane trees. However, the explicit formulas for a general $\eta$
are not expected to be as simple as those of two long cycles.

Next, in eq.~\eqref{eq:gen-new}, if we sum over all $\eta \vdash n$, it is not
hard to obtain
\begin{align}\label{cor:exc}
\sum_{\eta\vdash n}  \sum_{a\geq 0} a p^{\eta}_{\Lambda,a} =\bigl( n-\ell(\Lambda)\bigr) (n-1)! z_{\Lambda} - \sum_{i>0, \; j>0 \atop \Upsilon \rhd_{i, 2j+1} \Lambda} \kappa_{\Upsilon,\Lambda} (n-1)! z_{\Upsilon} \; .
\end{align}
We then realize that the left-hand side of eq.~\eqref{cor:exc} is the total number of exceedances of
$\mathfrak{p} \in \bigcup_{\eta \vdash n} U^{\eta}_{\Lambda}$ which can be easily obtained.

\begin{lemma}\label{lem:exc}
The total number of exceedances
\begin{align}\label{eq:exc}
\sum_{\eta \vdash n}  \sum_{a\geq 0} a p^{\eta}_{\Lambda, a} =\frac{n-m_1(\Lambda)}{2} (n-1)! z_{\Lambda} \; .
\end{align}
\end{lemma}

Combining eq.~\eqref{cor:exc} and eq.~\eqref{eq:exc}, and using $n+1-m_1(\Lambda)=\sum_{i>0} (i+1) m_{i+1}(\Lambda)$, we obtain the following proposition.

\begin{proposition} Let $\alpha \models n+1$. For any $\Lambda \vdash \alpha$, we have
\begin{align}\label{eq:baserecur}
\bigl( n+1-\ell(\Lambda) \bigr) z_{\Lambda} 
=\sum_{i>0, \; j>0 \atop \Upsilon \rhd_{i, 2j+1} \Lambda} \kappa_{\Upsilon,\Lambda} z_{\Upsilon} + \frac{ z_{\Lambda} }{2} \sum_{i>0} (i+1) m_{i+1} (\Lambda) \; .
\end{align}
\end{proposition}

In order to proceed, we need a few more notations.
Let $\lambda=1^{m_1} 2^{m_2} \cdots n^{m_n}\vdash n+1$.
For $i>0$ (and $m_{i+1}\neq 0$), denote $\lambda^{\downarrow (i+1)}$
the partition $\mu=1^{m_1} \cdots i^{m_i+1} (i+1)^{m_{i+1}-1} \cdots n^{m_n} \vdash n$, i.e.,~changing
an $i+1$ part to an $i$ part. 
Denote $\Lambda_i^{\downarrow (j+1)}$ the sequence of integer partitions
$(\lambda^{[1]}, \ldots, (\lambda^{[i]})^{\downarrow (j+1)}, \ldots, \lambda^{[k]})$.

As a consequence of eq.~\eqref{eq:long-new}, we have the following corollary.

\begin{corollary} Suppose $\alpha \models n+1$ and $\Lambda \vdash \alpha$. For  any $1\leq i \leq k$ and $j>0$, denote
$\Lambda_{i,j}=\frac{\alpha_i}{2} j m_j \bigl((\lambda^{[i]})^{\downarrow (j+1)}\bigr)$. Then, we obtain
\begin{multline}\label{eq:downarrow}
\bigl( n+1-\ell(\Lambda) \bigr) \Lambda_{i,j}  p^{(n)}_{\Lambda_i^{\downarrow (j+1)}} \\
=\sum_{t>0, \: d>0, \atop \Upsilon \rhd_{t,2d+1} \Lambda_i^{\downarrow (j+1)}} \kappa_{\Upsilon,\Lambda_i^{\downarrow (j+1)}}  \Lambda_{i,j}  p^{(n)}_{\Upsilon} 
 +\frac{(j+1) m_{j+1} (\lambda^{[i]})}{2} (n-1)! z_{\Lambda} \; .
\end{multline}
\end{corollary}

The upcoming lemma is a key ingredient enabling us to proceed further.

\begin{lemma} Let $\alpha \models n+1$ and $\Lambda \vdash \alpha$.
The following is true:
\begin{align}\label{eq:keylem}
\sum_{i>0, \; j>0} \Lambda_{i,j} \sum_{t>0,\; d>0 \atop \Upsilon \rhd_{t, 2d+1} \Lambda_i^{\downarrow (j+1)}} \kappa_{\Upsilon,\Lambda_i^{\downarrow (j+1)}}    p^{(n)}_{\Upsilon} 
= \sum_{t>0, \; d>0 \atop \Upsilon \rhd_{t, 2d+1} \Lambda} \kappa_{\Upsilon,\Lambda}  \sum_{i>0, \; j>0} \Upsilon_{i,j} p^{(n)}_{\Upsilon_i^{\downarrow (j+1)}} \; .
\end{align}
\end{lemma}

For any $\alpha \models n+1$ and $\Lambda \vdash \alpha$, we denote $\mathcal{T}_{\Lambda}= \sum_{i>0, \; j>0} \Lambda_{i,j}  p^{(n)}_{\Lambda_i^{\downarrow (j+1)}} $. Summing over all possible $i,j$ in eq.~\eqref{eq:downarrow} and applying eq.~\eqref{eq:keylem},
we obtain:
\begin{align}\label{eq:recur}
\bigl(n+1-\ell(\lambda) \bigr) \mathcal{T}_{\Lambda}=\sum_{\Upsilon \rhd_{t, 2d+1} \Lambda} \kappa_{\Upsilon, \Lambda} \mathcal{T}_{\Upsilon}
+\frac{(n-1)! z_{\Lambda}}{2} \sum_{i>0} (i+1) m_{i+1} (\Lambda) \; .
\end{align}
From now on, we also use the notation $\frac{\partial \alpha}{\partial i}:= (\alpha_1, \ldots, \alpha_i-1, \ldots, \alpha_k)$.
Now we are ready to prove our first main theorem.

\begin{theorem}\label{thm:main}
Let $\alpha \models n+1$ and $\Lambda \vdash \alpha$.
Suppose $\ell(\Lambda)$ has the same parity as $n$.
Then, we have
\begin{align}\label{eq:thm-main}
\mathcal{T}_{\Lambda} = (n-1)! z_{\Lambda} \; .
\end{align}
\end{theorem}
\proof Inspecting eq.~\eqref{eq:baserecur} and eq.~\eqref{eq:recur}, we observe that both sides of the equality
in the theorem satisfy the same recurrence.
Then, it suffices to compare the respective initial conditions. \qed

Let $\alpha\models n$ and $d=(d_1, d_2, \ldots, d_k)$ ($d_i>0$).
We denote $p^{(n)}_{\alpha, d}$ the number of pairs of $n$-cycles
whose product is $\alpha$-separated and the elements in $B_i$ are contained in $d_i$ cycles,
and we define $p^{(n)}_{\alpha}:= \sum_d p^{(n)}_{\alpha, d}$, i.e.,~the total number of $\alpha$-separating pairs.
As an application of Theorem~\ref{thm:main}, we have

\begin{theorem}\label{cor:main}
Let $\beta \models n+1$ and $d=(d_1, d_2, \ldots, d_k)$ where $\sum_i d_i$ has the same parity as $n$.
Then, we have
\begin{align}
\sum_{i=1}^k \frac{\beta_i (\beta_i-1)}{2} p^{(n)}_{\frac{\partial \beta}{\partial i} ,d } &= (n-1)! \prod_{i=1}^k C(\beta_i, d_i) \; ,\label{eq:cor-main1}\\
\sum_{i=1}^k \frac{\beta_i (\beta_i-1)}{2} p^{(n)}_{\frac{\partial \beta}{\partial i}  } &=  \frac{(n-1)! }{2} \prod_{i=1}^k \beta_i!\; . \label{eq:cor-main2}
\end{align}
\end{theorem}

Thanks to the above theorems, we obtain the results presented in Introduction.
Here is another main result.
For $\alpha\models n$ and $d=(d_1, d_2, \ldots, d_k)$,
we denote 
$$
Y_{d}^{\alpha}= (n-1)! C(\alpha_1+1, d_1)\prod_{i=2}^k C(\alpha_i, d_i).
$$
For $\alpha\models n$, we introduce the notation 
$\alpha^{(j)}=(\alpha_1+1, \ldots, \alpha_j-1, \ldots, \alpha_k) \models n$
for $j>1$. In addition, we inductively define $\alpha^{(j_1,\ldots, j_l)}:=\bigl( \alpha^{(j_1, \ldots, j_{l-1})}\bigr)^{(j_l)}$, for $j_t>1$.

\begin{theorem}\label{thm:explicit-gen}
For any $\alpha=(\alpha_1, \ldots, \alpha_k)\models n$ and $d=(d_1,\ldots, d_k)$,
we have
\begin{align}\label{eq:explicit-gen}
p^{(n)}_{\alpha, d}= \frac{Y^{\alpha}_d}{{\alpha_1+1 \choose 2}}+ \sum_{l>0} (-1)^l \sum_{(j_1,\ldots, j_l),\atop j_t>1}  \frac{Y^{\alpha^{(j_1,\ldots, j_l)}}_d }{{\alpha_1+l+1 \choose 2}} \prod_{t=1}^l   \left. {{1+\alpha^{(j_1,\ldots, j_t)}_{j_t}  \choose 2}} \right /{{\alpha_1+t \choose 2}}  .
\end{align}
\end{theorem}

The above results can be used to derive formulas concerning strong separation probabilities as well.

\section*{Acknowledgments}

I would like to thank Christian Reidys for encouragements and support.


\end{document}